\documentclass{amsart}

\usepackage{amssymb}
\usepackage{xcolor}
\usepackage{amsmath}
\usepackage{bm}

\newcommand{\half}{\frac{1}{2}}
\newcommand{\thalf}{\tfrac{1}{2}}
\newcommand{\summ}{\mathop{{\sum}^{\star}}}

\newcommand{\intt}{\int_{-\infty}^{\infty}}
\newcommand{\intx}{\int_{0}^{\infty}}

\makeatletter
\@namedef{subjclassname@2020}{\textup{2020} Mathematics Subject Classification}
\makeatother

\numberwithin{equation}{section}

\newtheorem{theorem}{Theorem}

\newtheorem{corollary}[theorem]{Corollary}

\begin{document}

\title{A reciprocity relation for the twisted second moment of the Riemman Zeta function}

\author{Rizwanur Khan}
\address{
Department of Mathematical Sciences\\ University of Texas at Dallas\\ Richardson, TX 75080-3021}
\email{rizwanur.khan@utdallas.edu }

\thanks{The author was supported by the National Science Foundation grants DMS-2344044 and DMS-2341239.}

\begin{abstract} We prove a reciprocity relation for the twisted second moment of the Riemann Zeta function. This provides an analogue to a formula of Conrey for Dirichlet $L$-functions.
\end{abstract} 
\maketitle

\section{Introduction}

Gauss's famous Law of Quadratic Reciprocity relates the Legendre symbol $(\frac{p}{q})$ to $(\frac{q}{p})$. Somewhat reminiscent of this, a reciprocity relation for $L$-functions is a relation connecting a (twisted) moment of $L$-functions (on the central line or at the central point) to a dual moment, where certain parameters are `inverted' on both sides. In an unpublished 2007 preprint, Conrey \cite{con} proved an elegant reciprocity relation for the twisted second moment of Dirichlet $L$-functions at the central point. This result garnered a lot of interest and was later extended by Young \cite{you}, Bettin \cite{bet}, Nordentoft \cite{nor}, and Djankovi\'{c} \cite{dja}. The relation, with one twist, is of the following shape. For primes $p<q$, we have
\begin{align}
\label{dirrec} \frac{1}{q^\half} \summ_{\chi \bmod q}  \chi(p)  |L(\thalf, \chi)|^2 = \text{main term } + \frac{1}{p^\half} \summ_{\chi \bmod p}  \chi(-q)  |L(\thalf, \chi)|^2+  \text{error term},
\end{align}
where the sums are restricted to primitive characters, the main term is roughly of size $(\frac{q}{p})^\half \log(\frac{q}{p})$ and the error term is roughly of size $O(\frac{p}{q})$. Thus the roles of $p$ and $q$ are inverted on both sides. The reciprocity relation \eqref{dirrec} is not an exact formula, unlike the so called spectral reciprocity relations referenced in \cite{humkha}. Although Bettin \cite[Theorem 2]{bet} did extend \eqref{dirrec} to a remarkable exact formula, it involved $L$-values off the central point.

In this paper we consider the twisted second moment of the Riemann Zeta function. Such a device has been extensively studied and has applications to our understanding of the zeros and size of the Riemann Zeta function. Some of the pioneering works in this regard are \cite{con2}, \cite{conghogon}, \cite{sou}. Asymptotics for the twisted second moment of the Riemann Zeta function (on average over the twists) have been given by Balasubramanian, Conrey, Heath-Brown \cite{balconhea} and Bettin, Chandee, Radziwi\l\l \ \cite{betcharad}. Our goal in this paper is to develop a {\it reciprocity relation} for the twisted second moment of the Riemann Zeta function. This sheds some light on the structure of the twisted second moment and provides an analogue to Conrey's relation for Dirichlet $L$-functions. We also mention that Bettin and Conrey \cite{betcon} proved a reciprocity relation for a certain sum involving cotangents which is related to the twisted second moment of the Riemann Zeta function.

Our twisted moment and dual moment do not have as symmetric of a relation to each other as \eqref{dirrec} does. However, this is to be expected. We recall the following spectral reciprocity identity, which can be obtained by the method in \cite{blokha}:
\begin{align*}
\sum_{f\in \mathcal{B}_0(1)} \lambda_f(p) \frac{L(\half, f)^4}{L(1, \text{ad} f)} {\phi}(t_f)+\ldots= \text{main term} +
\frac{1}{p^\half} \sum_{f\in \mathcal{B}_0(p)}\frac{L(\half, f)^4}{L(1, \text{ad} f)}\tilde{\phi}(t_f)+\ldots.
\end{align*}
where $\phi(t)$ is a nice weight function with a transform $\tilde{\phi}$, $\mathcal{B}_0(N)$ is an orthonormal basis of cuspidal Maass newforms $f$ on $\Gamma_0(N) \backslash\mathbb{H}$ with spectral parameter $t_f$, $\lambda_f(n)$ is the $n^\text{th}$ Hecke eigenvalue, and the ellipsis denotes analogous contributions from the holomorphic and Eisenstein series spectrum. The  family on the left hand side has an archimedean parametrization, but as a result of the twisting, the dual moment family has an arithmetic aspect. Following the same pattern, we might expect a twisted second moment of the Riemann Zeta function to be associated with a dual twisted second moment of Dirichlet $L$-functions, and this is what we find.
\begin{theorem}\label{main}
Let $p$ and $q$ be distinct odd primes. Let $T>1$. We have
\begin{align}
\label{thmline} \intt \Big(\frac{p}{q}\Big)^{it} |\zeta(\thalf+it)|^2 & \exp\Big(-\frac{t^2}{T^2}\Big)   dt  \\
\nonumber &=  \Big(\frac{\pi}{pq}\Big)^\half T\Big( \log \Big(\frac{T}{2\pi pq}\Big)+2\gamma+\half \frac{\Gamma'}{\Gamma}\Big(\half\Big)\Big)\\
\nonumber &\ +  \Big( \frac{T}{2\pi}\Big)^\half    \frac{p^\half }{p-1} \summ_{\substack{\chi\bmod p\\ \chi(-1)=1}}  \chi(q)  \ \intt \Gamma\Big(\frac{1-2it}{4}\Big) \Big( \frac{T}{2 q}\Big)^{it}  |L(\thalf+it,\chi)|^2 \ dt\\
\nonumber &\ +O \left((pqT)^\epsilon \Big(\frac{q}{p} \Big)^\half +  (pqT)^\epsilon \Big(\frac{p}{q} \Big)^\half  \right).
\end{align}
\end{theorem}
\noindent Since $\Gamma(\frac{1-2it}{4}) $ decays exponentially fast as $|t|\to\infty$, we should think of the dual moment integral as essentially being of finite length. We can also see that the dual moment is real, by pairing the contributions of $t$ and $-t$ and using that the complex conjugate of $\Gamma(\frac{1-2it}{4}) ( \frac{T}{2 q})^{it}$ is $\Gamma(\frac{1+2it}{4}) ( \frac{T}{2 q})^{-it}$, and pairing the contributions of $\chi$ and $\overline{\chi}$ when $\chi$ is not real.

Observe that the moment on the left hand side of \eqref{thmline} remains invariant if we interchange $p$ and $q$, by the substitution $t\to -t$ and the fact that $ \exp(-\frac{t^2}{T^2})$ is an even function. All other terms on the right hand side are also clearly symmetric in $p\leftrightarrow q$, except, at first glance, the dual moment. However we are reassured by the reciprocity relation \eqref{dirrec} which tells us that such expressions may be symmetric in $p\leftrightarrow q$. In fact, this {\it must} be the case as a direct consequence of our reciprocity relation \eqref{thmline}. So our reciprocity relation contains another one for free.
\begin{corollary}
Let $p$ and $q$ be distinct odd primes. Let $T>1$. We have
\begin{multline}
\label{corline}  \frac{p^\half }{p-1} \summ_{\substack{\chi\bmod p\\ \chi(-1)=1}}  \chi(q)  \ \intt \Gamma\Big(\frac{1-2it}{4}\Big) \Big( \frac{T}{2 q}\Big)^{it}  |L(\thalf+it,\chi)|^2 \ dt \\
 =  \frac{q^\half }{q-1} \summ_{\substack{\chi\bmod q\\ \chi(-1)-1}}  \chi(p)  \ \intt  \Gamma\Big(\frac{1-2it}{4}\Big) \Big( \frac{T}{2 p}\Big)^{it}  |L(\thalf+it,\chi)|^2 \ dt+O \left((pqT)^\epsilon \Big(\frac{q}{pT} \Big)^\half +  (pqT)^\epsilon \Big(\frac{p}{qT} \Big)^\half  \right).
\end{multline}
\end{corollary}
\noindent Unlike \eqref{dirrec}, we do not have a main term because of the extra oscillation arising from the integral. Also note that the error term $O((pqT)^\epsilon (\frac{p}{qT})^\half )$ is the expected order of magnitude for the moment on the left hand side of \eqref{corline}, while $O((pqT)^\epsilon (\frac{q}{pT})^\half )$ is the expected order of magnitude for the dual moment on the right hand side. Such expectations however are known in only limited ranges, as we discuss in the next paragraph, so \eqref{corline} is non-trivial.

Going back to the reciprocity relation \eqref{dirrec}, one interesting consequence is the following. While an asymptotic is known \cite[Corollary 2]{bet} for the moment on the left hand side for $p<q^{\half-\epsilon}$, by \eqref{dirrec} we have an asymptotic for the {\it difference} between the moment on the left hand side and the dual moment on the right hand side in the wider range $p<q^{1-\epsilon}$. To see that we get the analogue of this for the Riemann Zeta function, suppose that $p\asymp q$. While an asymptotic is known for the moment on left hand side of \eqref{thmline} only for $p<T^{\half-\epsilon}$, we get an asymptotic for the {\it difference} between the moment on the left hand side and the dual moment on the right hand side in the wider range $p<T^{1-\epsilon}$. 

We have tried to present the reciprocity relation in a straightforward form, and therefore opted to make some simplifications. We have worked with prime twists in Theorem \ref{main}, although this restriction is not a serious one. For general integer value twists, the entire proof would go through in the same way but the orthogonality property of primitive Dirichlet characters, used in \eqref{dirchar} and \eqref{2nd}, would be more complicated than the prime case. We have also worked with the specific weight function  $\exp(-\frac{t^2}{T^2})$ on the left hand side of \eqref{thmline} so that all transforms which arise can be explicitly evaluated. 

Instead of considering a second moment of the Riemann Zeta function along a dyadic interval with twists as long as possible, one may consider the opposite scenario of an untwisted second moment along as short an interval as possible, say using a weight function $\exp(-\frac{(t-T)^2}{H^2})$ with $H\le T$. The dual moment in this case would be a second moment of the Riemann Zeta function against a transform function with (inverted) support of size $\frac{T}{H}$. We do not pursue this as it may already be derivable from an Atkinson type formula \cite[Theorem 4.1]{mot}.

\

{\bf Acknowledgement.} I thank Matt Young and Joseph Leung for some illuminating discussions.

\section{Sketch}

To illustrate our idea for Theorem \ref{main}, we sketch a new proof of \eqref{dirrec} which is amenable to generalization to give an analogous relation for the Riemann Zeta function. Start with
\[
 \summ_{\chi \bmod q}  \chi(p)  |L(\thalf, \chi)|^2
\]
and apply the functional equation to one of the central values, to get something involving
\[
 \summ_{\chi \bmod q}  \chi(p) \tau(\chi) L(\thalf, \overline{\chi})^2 =  \summ_{\chi \bmod q}  \chi(p) \tau(\chi) \sum_{n,m\ge 1} \frac{\overline{\chi}(nm)}{(nm)^\half},
 \]
after expanding $L(\thalf, \overline{\chi})^2$ as a Dirichlet series, where $\tau(\chi)$ is the Gauss sum. Exchanging the order of summation, opening the Gauss sum by writing
 \[
 \tau(\chi)=\sum_{a\bmod q} e\Big(\frac{a}{q}\Big) \chi(a),
 \]
executing the sum over $\chi$ and sum over $a$, leads to something involving
 \[
\sum_{n,m\ge 1} e\Big(\frac{nm\overline{p}}{q}\Big)\frac{1}{\sqrt{nm}}.
 \]
 Similarly, if we start with 
 \[
 \summ_{\chi \bmod p}  \chi(-q)  |L(\thalf, \chi)|^2 ,
\]
we would end up with
 \[
\sum_{n,m\ge 1} e\Big(\frac{-nm\overline{q}}{p}\Big)\frac{1}{\sqrt{nm}},
 \]
 where $\overline{q}$ denotes the multiplicative inverse of $q$ modulo $p$. 
 The two expressions are `approximately' equal, using the additive reciprocity identity (which follows from the Chinese Remainder Theorem),
 \[
 e\Big(\frac{nm\overline{p}}{q}\Big) = e\Big(\frac{-nm\overline{q}}{p}\Big)e\Big(\frac{nm}{pq}\Big).
 \]
 We can take the Taylor series expansion of $e(\frac{nm}{pq})$, and use only the first term to get the reciprocity formula with an error term. Thus the steps are: applying the functional equation, opening the oscillatory factor from the functional equation, averaging over the family, an approximation step, and then repeating the preceding steps in reverse.

\section{Proof of Theorem \ref{main}}

\subsection{Set up}
Define, for $\Re(s)>-\half$, where $\Re$ denotes the real part,
\begin{align*}
F(s)= \intt \Big(\frac{p}{q}\Big)^{it} \zeta(\thalf-s-it)\Big(\zeta(\thalf+s+it) - (s-\thalf+it)^{-1} \Big) \exp\Big(-\frac{t^2}{T^2}\Big) dt,
\end{align*}
and note that this is analytic in the specified domain. Note that $F(0)$ is essentially the moment on the left hand side of \eqref{thmline}. More precisely, we have
\[
F(0)= \intt \Big(\frac{p}{q}\Big)^{it} |\zeta(\thalf+it)|^2 \exp\Big(-\frac{t^2}{T^2}\Big)  dt + O(T^\epsilon).
\]
The plan now is to obtain a different expression for $F(s)$ through analytic continuation, then to set $s=0$ in that new expression.
\subsection{Functional equation} By the functional equation,
\begin{align*}
F(s)=\intt \Big(\frac{p}{q}\Big)^{it}  \pi^{-s- i t} \frac{ \Gamma( \frac14+\frac{s}{2}+\frac{it}{2}) }{ \Gamma( \frac14-\frac{s}{2}-\frac{it}{2})} \zeta(\thalf+s+it) \Big(\zeta(\thalf+s+it) - (s-\thalf+it)^{-1} \Big) \exp\Big(-\frac{t^2}{T^2}\Big)  dt
\end{align*}
and by the reflection and duplication formulae of the gamma function, we have
\begin{multline*}
F(s)= 2 \intt  \Big(\frac{p}{q}\Big)^{it}  (2\pi)^{-\half-s- i t}  \cos(\pi(\tfrac14 + \tfrac{s}{2} +\tfrac{it}{2}))  \Gamma( \thalf +s +it) \zeta(\thalf+s+it) \Big(\zeta(\thalf+s+it) - (s-\thalf+it)^{-1} \Big) \\
\times  \exp\Big(-\frac{t^2}{T^2}\Big)  dt.
\end{multline*}
\subsection{Opening an oscillatory factor}  For $|\Re(s)|<\half$, we write $\cos(\pi(\tfrac14 + \tfrac{s}{2} +\tfrac{it}{2}))  \Gamma( \thalf +s +it) $ as a Mellin transform using \cite[17.43.3]{gr}  to get
\begin{multline}
\label{ready-interchange} F(s)= 2 \intt   \intx     x^{-\half+s+it} \cos x \ dx \   \Big(\frac{p}{q}\Big)^{it}  (2\pi)^{-\half-s- i t}   \zeta(\thalf+s+it) \Big(\zeta(\thalf+s+it) - (s-\thalf+it)^{-1} \Big) \\
\times \exp\Big(-\frac{t^2}{T^2}\Big)  dt.
\end{multline}
Integrating by parts and using $|\sin x|\le \min\{1,|x|\}$, we have
\begin{multline*}
F(s)=  -2 \intt  \intx     x^{-\frac32+s+it} \sin x \ dx \   \Big(\frac{p}{q}\Big)^{it}   (2\pi)^{-\half-s- i t}  (-\thalf+s+it) \zeta(\thalf+s+it) \\
\times \Big(\zeta(\thalf+s+it) - (s-\thalf+it)^{-1} \Big) 
 \exp\Big(-\frac{t^2}{T^2}\Big)  dt.
\end{multline*}
This $x$-integral conditionally converges for $-\half < \Re(s) < \frac32$. 
Integrating by parts, we write
\[
 \intx     x^{-\frac32+s+it} \sin x \ dx = \int_0^{\frac{\pi}{2}}   x^{-\frac32+s+it} \sin x \ dx + \int_{\frac{\pi}{2}}^{\infty}   (-\tfrac32+s+it) x^{-\frac52+s+it} \cos x \ dx,
\]
where both integrals on the right hand side are absolutely convergent for $-\half < \Re(s) < \frac32$. Now we restrict to $\half <\Re (s)<\frac32$ and expand $\zeta(\thalf+s+it)$ into a absolutely convergent Dirichlet series. Thus we get, after exchanging the order of integration and summation, 
\begin{align*}
 F(s)= F_{1}(s)+F_{2}(s)+F_{3}(s),
\end{align*}
where
\begin{align}
\label{f1} F_1(s) = -2\sum_{n,m\ge 1} \frac{s-\half}{(2\pi nm)^{\half+s}}  \Bigg( \int_0^{\frac{\pi}{2}}     I_1(x)   \sin x  \ dx
+ \int_{\frac{\pi}{2}}^\infty     I_2(x)  \cos x \ dx \Bigg),
\end{align}
for
\begin{align*}
&I_1(x)=  \intt x^{-\frac32+s}   \Big( \frac{xp}{2\pi nmq} \Big)^{it}     \exp\Big(-\frac{t^2}{T^2}\Big)   dt,\\
&I_2(x) = \intt (-\tfrac32+s+it) x^{-\frac52+s}   \Big( \frac{xp}{2\pi nmq} \Big)^{it}     \exp\Big(-\frac{t^2}{T^2}\Big)   dt,
\end{align*}
and $F_2(s)$ and $F_3(s)$ are to be defined presently. We integrate by parts the second $x$-integral in \eqref{f1}, by differentiating $\cos x$ and integrating $I_2(x)$ (under the $t$-integral sign). This gives
\begin{align*}
F_1(s)&= -2\sum_{n,m\ge 1} \frac{s-\half}{(2\pi nm)^{\half+s}}   \int_0^{\infty}     I_1(x)   \sin x  \ dx\\
&= -2\sum_{n,m\ge 1} \frac{s-\half}{(2\pi nm)^{\half+s}}   \int_0^{\infty}   x^{-\frac32+s}   \sin x    \intt  \Big( \frac{xp}{2\pi nmq} \Big)^{it}     \exp\Big(-\frac{t^2}{T^2}\Big)   dt \ dx
\end{align*}
The integrals and sum all converge absolutely, as we will see. 

At this point we return to the definitions of $F_1(s)$ and $F_3(s)$ before proceeding further. By the same integration by parts step, we have
\begin{align}
\label{f2} F_2(s) =  &-2\sum_{n,m\ge 1} \frac{1}{(2\pi nm)^{\half+s}}   \int_0^{\infty}   x^{-\frac32+s}   \sin x    \intt  it \Big( \frac{xp}{2\pi nmq} \Big)^{it}     \exp\Big(-\frac{t^2}{T^2}\Big)   dt \ dx,\\
\nonumber F_3(s)=&2\sum_{n,m\ge 1} \frac{1}{(2\pi n)^{\half+s}}   \int_0^{\infty}   x^{-\frac32+s}   \sin x    \intt  \Big( \frac{xp}{2\pi nq} \Big)^{it}     \exp\Big(-\frac{t^2}{T^2}\Big)   dt \ dx.
\end{align}

\subsection{Averaging over the family (evaluating the integral on the critical line)}
Using the Fourier transform \cite[17.23.13]{gr}
\begin{align*} 
 \pi^{-\half}\intt  e^{it\xi}     \exp\Big(-\frac{t^2}{T^2}\Big)   dt = T \exp\Big( -\frac{T^2}{4}\xi^2  \Big)
 \end{align*}
with $\xi=\log(\frac{xp}{2\pi nmq})$, we get
\begin{align*}
F_1(s) = -2 \pi^\half T\sum_{n,m\ge 1} \frac{s-\half}{(2\pi nm)^{\half+s}}   \int_0^\infty   x^{-\frac32+s}   \sin x   \exp\left( -\frac{T^2}{4}\Big(\log \frac{xp}{2\pi nmq}\Big)^2 \right) \ dx,
\end{align*}
and note that the $x$-integral here converges absolutely. Making the substitution $y=\frac{T}{2}\log(\frac{xp}{2\pi nmq})$, we get
\begin{align*}
F_1(s) =  & -4 \pi^\half  \sum_{n,m\ge 1} \frac{s-\half}{(2\pi nm)^{\half+s}}  \intt \Big(\frac{2 \pi nmq}{p} e^{\frac{2y}{T}}\Big)^{-\frac12+s}   \sin \Big(\Big(\frac{2\pi nmq}{p}\Big) e^{\frac{2y}{T}}\Big)  e^{-y^2} \ dy\\
=  & -4 \pi^\half  \Big(\frac{q}{p} \Big)^{-\frac12+s}  \sum_{n,m\ge 1} \frac{s-\half}{2\pi nm}  \intt  e^{\frac{2y}{T}(-\frac12+s)}  \sin \Big(\Big(\frac{2\pi nmq}{p}\Big) e^{\frac{2y}{T}}\Big)  e^{-y^2} \ dy.
\end{align*}
The sum converges absolutely by integration by parts of the integral. In fact, up to an error of $O_\epsilon((pqT)^{-100})$, the sum may be restricted to $nm\le (pqT)^\epsilon \frac{Tp}{q}$. Thus we can analytically continue $F_1(s)$ to an entire function, and at $s=0$ we have
\begin{align*}
F_1(0) =   \frac{1}{\pi^\half}  \Big(\frac{q}{p} \Big)^{-\frac12} \sum_{n,m\ge 1} \frac{1}{ nm}  \intt e^{-\frac{y}{T}}  \sin \Big(\Big(\frac{2\pi nmq}{p}\Big) e^{\frac{2y}{T}}\Big)  e^{-y^2} \ dy,
\end{align*}
and the bound 
\[
|F_1(0)|\ll (pqT)^\epsilon \Big(\frac{p}{q}\Big)^\half.
\]
The treatment for $F_3(s)$ is similar: we can continue it to $s=0$ and we get the bound $F_3(0)\ll (pqT)^\epsilon (\frac{p}{q})^{\frac12}$.

Now we turn to $F_2(s)$. This time we insert
\begin{align}
\label{fouriertransform2} \pi^{-\half}\intt  it e^{it\xi} \exp\Big(-\frac{t^2}{T^2}\Big)   dt = -\frac{T^3}{2}  \xi \exp\Big( -\frac{T^2}{4} \xi^2 \Big)
\end{align}
with $\xi=\log( \frac{xp}{2\pi nmq})$ into the definition \eqref{f2}, and after the same substitution $y=\frac{T}{2}\log(\frac{xp}{2\pi nmq})$ and continuing to $s=0$, we get
\begin{align}
\nonumber F_2(0) &=    \frac{2T}{\pi^\half}  \Big(\frac{q}{p} \Big)^{-\frac12} \sum_{n,m\ge 1} \frac{1}{ nm}  \intt  e^{-\frac{y}{T}} \sin \Big(\Big(\frac{2\pi nmq}{p}\Big) e^{\frac{2y}{T}}\Big)  y e^{-y^2} \ dy\\
\label{f2(0)} &= \Im  \  \frac{2T}{\pi^\half}  \Big(\frac{q}{p} \Big)^{-\frac12} \sum_{n,m\ge 1} \frac{1}{ nm}  \intt    e\Big(\Big(\frac{ nmq}{p}\Big) e^{\frac{2y}{T}}\Big)  y e^{-y^2-\frac{y}{T}} \ dy,
\end{align}
where $e(x)=e^{2\pi i x}$.

\subsection{An approximation} Up to an error of $O_\epsilon((pqT)^{-100})$, by repeated integration by parts, the sum in \eqref{f2(0)} may be restricted to $nm\le (pqT)^\epsilon \frac{Tp}{q}$, and the integral may be restricted to $|y|\le (pqT)^\epsilon$. Then, we can make the Taylor approximations
\begin{align*}
&e^{-\frac{y}{T}}=1+O\Big(\frac{y}{T}\Big)=1+O\Big(\frac{(pqT)^\epsilon}{T}\Big) ,\\
&e\Big(\Big(\frac{ nmq}{p}\Big) e^{\frac{2y}{T}}\Big)=e\Big(\frac{ nmq}{p}\Big)e\Big(\frac{ nmq}{p}\frac{2y}{T}\Big)+O\Big(\frac{(pqT)^\epsilon}{T}\Big) ,
\end{align*}
to get
\begin{align}
\label{afterapprox} F_2(0) &= \Im  \  \frac{2T}{\pi^\half} \sum_{\substack{n,m\ge 1\\ nm\le (pqT)^\epsilon \frac{pT}{q}}} \frac{1}{ nm} e\Big(\frac{ nmq}{p}\Big) \int_{-(pqT)^\epsilon}^{(pqT)^\epsilon} \Big(\frac{q}{p} \Big)^{-\frac12}    e\Big(\frac{ nmq}{p}\frac{2y}{T}\Big)  y e^{-y^2} \ dy + O\Big((pqT)^\epsilon \Big(\frac{p}{q}\Big)^{\frac12}\Big).
\end{align}
The integral and sum may now be extended back to infinity, up to a negligible error, after which we can evaluate the integral (which completes the averaging over the family because this integral originated from opening the oscillatory factor) using \eqref{fouriertransform2} with $\xi= \frac{4\pi nmq}{pT}$. We get
\begin{align*}
F_2(0) &= \Re  \  4\pi    \Big(\frac{q}{p} \Big)^{\frac12}  \sum_{n,m\ge 1}  e\Big(\frac{ nmq}{p}\Big)  \exp\Big(-\Big(\frac{ 2\pi nm q}{pT}\Big)^2\Big) + O\Big((pqT)^\epsilon \Big(\frac{p}{q}\Big)^{\frac12}\Big)\\
&= 4\pi    \Big(\frac{q}{p} \Big)^{\frac12}  \sum_{n,m\ge 1}  \cos\Big(\frac{ 2\pi nmq}{p}\Big)  \exp\Big(-\Big(\frac{ 2\pi nm q}{pT}\Big)^2\Big) + O\Big((pqT)^\epsilon \Big(\frac{p}{q}\Big)^{\frac12}\Big).
\end{align*}

\

\

Now to get to the dual moment, we repeat the preceding steps, in reverse. Thus the titles of following subsections indicate the {\it reverse} of what is done.

\subsection{Averaging over the family (evaluating the character sum and integral on the critical line)}
Splitting the sum according to the residue class of $nm$ mod $p$, we have
\begin{align}
\nonumber F_2(0) = 4\pi    \Big(\frac{q}{p} \Big)^{\frac12}   \Bigg(&   \sum_{\substack{n,m\ge 1\\nm\equiv 0\bmod p}}    \exp\Big(-\Big(\frac{ 2\pi nm q}{pT}\Big)^2\Big) \\
\label{sums} & +   \summ_{r\bmod p} \cos \Big(\frac{2\pi rq}{p}\Big) \sum_{\substack{n,m\ge 1\\nm\equiv r\bmod p}}   \exp\Big(-\Big(\frac{ 2\pi nm q}{pT}\Big)^2\Big)\Bigg) + O\Big((pqT)^\epsilon \Big(\frac{p}{q}\Big)^{\frac12}\Big),
\end{align}
where $\summ$ indicates that the sum is over the non-zero residue classes. We rewrite the sum with the condition $nm\equiv 0\bmod p$ as a sum with $n=pn'$ plus a sum with $m=pm'$ minus a sum with $nm=p^2n'm'$. For the other sum in \eqref{sums}, we use the orthogonality of Dirichlet characters to pick out the condition $nm\equiv r\bmod p$. We get
\begin{align}
\nonumber F_2(0) = &4\pi   \Big(\frac{q}{p} \Big)^{\frac12}  \Bigg(   2  \sum_{\substack{n,m\ge 1}}    \exp\Big(-\Big(\frac{ 2\pi nm q}{T}\Big)^2\Big)   - \sum_{\substack{n,m\ge 1}}    \exp\Big(-\Big(\frac{ 2\pi nm pq}{T}\Big)^2\Big) \\
\label{dirchar}  & + \frac{1}{p-1} \summ_{r\bmod p} \cos \Big(\frac{2\pi rq}{p}\Big) \sum_{\chi\bmod p} \overline{\chi}(r) \sum_{\substack{n,m\ge 1}}  \chi(nm) \exp\Big(-\Big(\frac{ 2\pi nm q}{pT}\Big)^2\Big)\Bigg) + O\Big( (pqT)^\epsilon \Big(\frac{p}{q}\Big)^{\frac12}\Big).
\end{align}
Next, replacing $r$ with $r\overline{q}$, using the Mellin transform pair \cite[17.43.2]{gr}, writing
\[
\sum_{n,m\ge 1} \frac{1}{(nm)^w} = \zeta(w)^2, \sum_{n,m\ge 1} \frac{\chi(nm)}{(nm)^w} = L(w, \chi)^2, \sum_{n,m\ge 1} \frac{\chi_0(nm)}{(nm)^w} = \Big(1-\frac{1}{p^{w}}\Big)^2\zeta(w)^2,
\]
where $\Re(w)>1$ and $\chi_0$ is the trivial character, and using
$
\summ\limits_{r\bmod p} \cos (\frac{2\pi r}{p})=-1,
$
we get
\begin{align}
\nonumber F_2(0) =     &  \Big(\frac{q}{p} \Big)^{\frac12} \frac{1}{2\pi i} \int\limits_{(2)} 2\pi \Gamma\Big(\frac{w}{2}\Big) \Big( \frac{T}{2\pi}\Big)^w \Big(\frac{2}{q^w}- \frac{1}{(pq)^w} - \frac{1}{p-1}\frac{p^w}{q^w} \Big(1-\frac{1}{p^{w}}\Big)^2 \Big)\zeta(w)^2 \ dw\\
\label{2nd}  +   & \Big(\frac{q}{p} \Big)^{\frac12}  \frac{1}{2\pi i}  \int\limits_{(2)} 2\pi  \Gamma\Big(\frac{w}{2}\Big) \Big( \frac{Tp}{2\pi q}\Big)^w \frac{1}{p-1} \summ_{\chi\bmod p}   \chi(q) \summ_{r\bmod p} \cos \Big(\frac{2\pi r}{p}\Big) \overline{\chi}(r) L(w,\chi)^2  \ dw\\
\nonumber +& O\Big((pqT)^\epsilon \Big(\frac{p}{q}\Big)^{\frac12}\Big).
\end{align}
For the first integral, we move the line of integration to $\Re(w)=\epsilon$, crossing a double pole of $\zeta(w)^2$ at $w=1$. The residue from this pole is the main term displayed in \eqref{thmline} and the shifted integral is of size $O( (pqT)^\epsilon (\frac{q}{p})^{\frac12})$.

\subsection{Opening the oscillatory factor (Gauss sum)}

In the second integeral \eqref{2nd}, we have
\[
 \summ_{r\bmod p} \cos \Big(\frac{2\pi r}{p}\Big) \overline{\chi}(r) =0
\]
if $\chi(-1)=-1$, by pairing up the contributions of $r$ and $-r$ and using that cosine is an even function, and 
\[
 \summ_{r\bmod p} \cos \Big(\frac{2\pi r}{p}\Big) \overline{\chi}(r) =   \summ_{r\bmod p} \Big(\cos \Big(\frac{2\pi r}{p}\Big) + i \sin \Big(\frac{2\pi r}{p}\Big) \Big) \overline{\chi}(r)= \summ_{r\bmod p} e \Big(\frac{ r}{p}\Big) \overline{\chi}(r)  = \tau(\overline{\chi})
\]
if $\chi(-1)=1$, by pairing up the contributions of $r$ and $-r$ to see that $ \summ\limits_{r\bmod p} \sin(\frac{ 2\pi r}{p})  \overline{\chi}(r) $ vanishes.

\subsection{Applying the functional equation}
Moving the second integral \eqref{2nd} to the line $\Re(w)=\frac12$, writing $w=\half+it$ so that $dw=idt$, applying the functional equation of the Dirichlet $L$-functions, and using $|\tau(\chi)|^2=p$, we get that it equals
\begin{align*}
& \Big( \frac{T}{2\pi}\Big)^\half  \frac{1}{p-1} \summ_{\substack{\chi\bmod p\\\chi(-1)=1}}  \chi(q)  \  \intt \Gamma\Big(\frac{1+2it}{4}\Big) \Big( \frac{Tp}{2\pi q}\Big)^{it}   \tau(\overline{\chi}) L(\thalf+it,\chi)^2 \ dt\\
=&  \Big( \frac{T}{2\pi}\Big)^\half    \frac{p^\half }{p-1} \summ_{\substack{\chi\bmod p\\\chi(-1)=1}}  \chi(q)  \ \intt \Gamma\Big(\frac{1-2it}{4}\Big) \Big( \frac{T}{2 q}\Big)^{it}  |L(\thalf+it,\chi)|^2 \ dt.
\end{align*}

\

\

\bibliographystyle{amsplain} 
\bibliography{2nd-moment-zeta}

\end{document}